\documentclass[11pt]{article}
\usepackage{authblk}
\usepackage{graphicx}
\usepackage{amssymb}
\usepackage{amsmath}
\usepackage{amsthm}
\usepackage{tikz}
\usepackage{booktabs}
\usepackage{array}
\usepackage{multirow}
\usepackage{float}
\usepackage[round]{natbib}
\usepackage{xurl}
\usepackage{hyperref}
\usepackage{titlesec}
\usepackage{etoolbox}

\appto\appendix{
    \titleformat{\section}[block]
    {\normalfont\Large\bfseries}
    {APPENDIX \thesection:}
    {0.5em}
    {}
}

\usetikzlibrary{matrix}

\newtheorem{definition}{Definition}
\newtheorem{theorem}{Theorem}
\newtheorem{corollary}{Corollary}
\newtheorem{lemma}{Lemma}
\newtheorem{example}{Example}
\newtheorem{proposition}{Proposition}

\begin{document}

\title{%\vspace*{2cm}
Hierarchical Traffic Networks: A Framework for Reducing Braess' Paradox}
\author[1]{Guilherme Goncalves Moreira}
\author[2]{Vadim Zverovich
\footnote{Corresponding author (e-mail: vadim.zverovich@uwe.ac.uk)}
}
\affil[ ]{\footnotesize $^1$Warwick Mathematics Institute, University of Warwick,
Coventry, UK}
\affil[ ]{\footnotesize $^2$Mathematics and Statistics Research Group, University of the West of England, Bristol, UK}
\date{12 November 2026}

\maketitle

\begin{abstract}
Braess' paradox occurs when adding a link to a traffic network increases travel times at user equilibrium. We investigate its occurrence in the classical Braess network whose linear travel-time functions of links incorporate a hierarchical parametrisation. 
We derive necessary and sufficient conditions for the paradox to occur and also sufficient conditions when one hierarchy-scaling parameter is sufficiently large. For fixed free-flow scaling, we bound the probability that randomly assigned hierarchy parameters admit the paradox for some remaining scaling parameter and total flow. Finally, we use Monte Carlo simulations to estimate the paradox's probability when scaling parameters and hierarchy levels are sampled from specified distributions. The estimated probabilities are generally very low and depend on the sampling distributions. Compared with non-hierarchical networks, these lower probabilities suggest that hierarchical networks may better model real-world transportation systems from a congestion-minimisation perspective.
\end{abstract}

\section{Introduction}

The mathematical modelling of traffic systems developed substantially during the 1950s. Notably, \citet{Wardrop1952} introduced the principle of user equilibrium, according to which traffic distributes itself across alternative routes such that no traveller can reduce their travel time by unilaterally changing routes. \citet{Beckmann1956} subsequently formulated the user-equilibrium problem as a non-linear optimisation problem, providing a mathematical framework for computing equilibrium traffic flows.

A particularly surprising consequence of user equilibrium behaviour was later discovered. Braess' paradox was first described by the German mathematician Dietrich Braess in a paper originally published in 1968 \citep[English translation]{DBraess2005}. Braess demonstrated that adding a new link to a transport network can increase users' travel times at equilibrium
when users independently choose routes that minimise their own travel times. \citet{JDMurchland1970} subsequently discussed the paradox and reported an early real-world example consistent with it. Murchland described how major road investments in Stuttgart city centre failed to produce the expected improvements until a cross-street was closed to traffic. 

Subsequent work sought to characterise the conditions under which the paradox occurs in the classical network introduced by Braess. \citet{PasPrincipio1997} showed that its occurrence depends on the parameters of the link congestion functions and the total travel demand (i.e. total flow). For a network with fixed congestion functions, they demonstrated that the paradox can occur only when demand falls within a certain intermediate range. However, their analysis assumed identical travel-time functions for two pairs of links, a symmetry condition that is uncommon in real transport networks. \citet{ZverovichAvineri2015} relaxed this assumption by considering arbitrary linear travel-time functions and formulating conditions for the occurrence of the paradox in a generalised network configuration.

Research also broadened the study of Braess' paradox by altering the underlying network structure, behavioural assumptions and definition of user cost. \citet{Yang1998Capacity} introduced explicit link capacities into the classical network and identified a ‘capacity paradox’, in which adding a new link can reduce the maximum demand that the network can accommodate at user equilibrium. \citet{ValiantRoughgarden2006} proved that Braess' paradox is likely to occur in a natural random-network model. More precisely, they showed that, with high probability as the number of vertices tends to infinity, a traffic demand and a set of edges can be identified such that removing those edges reduces equilibrium latency by a constant factor. Their analysis also showed that equilibrium behaviour in their
large random-network model resembles that of Braess' original network. This supports the relevance of Braess' original network configuration as a theoretical benchmark. 

\citet{XuanXiaoHenry2014BRUE} relaxed the assumption of perfect rationality by allowing travellers to choose any route whose travel cost lies within an ‘indifference band’ of the shortest-path cost. They derived conditions relating the occurrence of the paradox to the width of this indifference band and the level of travel demand. More recently, \citet{Cheng2020Pollution} incorporated both the value of travel time and travellers’ perceived health damage from pollution into route choice, deriving conditions under which a related transportation network paradox occurs. Collectively, these studies demonstrate that such counterintuitive network effects are not confined to the assumptions or structure of Braess' original model.

%Beyond these extensions, other studies have examined the detection, analysis, and practical implications of Braess' paradox in transportation networks. \citet{Bagloee2019} propose a method for detecting the paradox in the context of no-car zone policies, while \citet{Bazzan2005} use simulation to study route recommendation and learning in abstract and microscopic models. \citet{Zhao2014} explore the paradox under stochastic user equilibrium and its implications for network robustness. Earlier contributions include \citet{Fisk1981}, who provide empirical evidence for equilibrium paradoxes, and \citet{Penchina1997}, who derive the maximum penalty in a minimal critical network.

Beyond these extensions, other studies have examined the detection, analysis, and practical implications of Braess' paradox in transportation networks. \citet{Bagloee2019} investigate the paradox in the context of no-car zone policies, proposing a method for detecting Braess' paradox and assessing whether such policies are beneficial. \citet{Bazzan2005} examine the paradox through simulation, comparing route recommendation and learning behaviour in both abstract and microscopic traffic models. \citet{Zhao2014} analyse the paradox under stochastic user equilibrium, highlighting the relationship between Braess' paradox and the robustness of traffic networks. Earlier work by \citet{Fisk1981} provides empirical evidence for equilibrium paradoxes and discusses their implications for optimal planning strategies. Finally, \citet{Penchina1997} identifies the maximum penalty that can arise in a minimal critical network, offering a quantitative bound on the worst-case impact of the paradox.

Braess-type effects have also been investigated in specialised transportation settings and other networked systems. Within transportation, \citet{Zhang2024} examined the occurrence of the paradox in autonomous-vehicle networks when parking-related travel is considered. In other networked systems, \citet{Nagurney2007Internet} formulated a time-dependent version of the paradox in an Internet model with multiple classes of traffic, while \citet{Schafer2022} experimentally demonstrated the paradox in a laboratory-scale AC power grid and showed through detailed simulations that planned extensions to the German high-voltage grid could induce the same effect.

Most directly relevant to the present study, Chapter~2 of \citet{VZverovich2021} developed a probabilistic framework for analysing Braess' paradox in the classical network configuration. By treating the parameters of the travel-time functions as random variables, the study examined the probability of the paradox under many distributions, with particular emphasis on the Erlang distribution. 
The present paper extends this probabilistic approach to network configurations incorporating the hierarchical parametrisation of \citet{Avineri2016} and investigates how the distributions of the network parameters affect the probability of the paradox occurring.

Although Braess' paradox has been studied extensively in classical network configurations, real road networks exhibit more complex and hierarchical structures.
Such hierarchies can emerge through interactions between local investment rules, traveller behaviour and the underlying network topology \citep{Yerra2005Hierarchies}. For example, roads in Great Britain are organised through a hierarchical classification system of motorways, A roads, B roads and local roads, reflecting differences in their intended functions within the network \citep{DfT2012RoadClassification}. This raises the question of whether the paradox persists in hierarchical road networks and how likely it is to occur under different assumptions about the network parameters.

To address these questions, we investigate Braess' paradox in the classical network configurations incorporating the hierarchical parametrisation of \citet{Avineri2016}. Hierarchy is represented through the link travel-time parameters, while the topology of Braess' network is retained. 
In Section~\ref{SectionHN},
we derive necessary and sufficient conditions for the paradox to occur and describe some situations when the paradox cannot occur.
Section~\ref{Section4} provides sufficient conditions for the paradox to occur when one hierarchy-scaling parameter is sufficiently large. For fixed free-flow scaling, we also derive bounds on the probability that randomly assigned hierarchy parameters admit the paradox for some choice of the remaining scaling parameter and total flow. The section concludes with two explicit examples of Braess' paradox in a hierarchical network.
In Section~\ref{Section5},
we use Monte Carlo simulations to estimate the paradox's probability when scaling parameters and hierarchy levels are sampled from specified distributions. The estimated probabilities are generally very low and depend on the sampling distributions. Compared with non-hierarchical networks, these lower probabilities suggest that hierarchical networks may better model real-world transportation systems from a congestion-minimisation perspective.

\section{Classical setup}

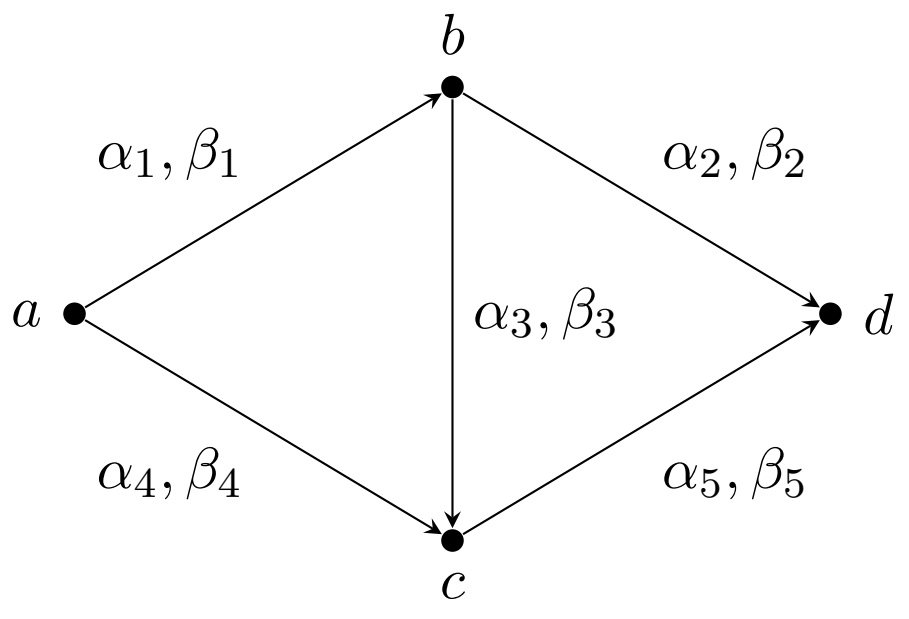
\begin{figure}[t]
\centering
\begin{tikzpicture}[
  >=stealth,
  node distance=3cm,
  vertex/.style={circle, fill=black, inner sep=1.5pt, label distance=3pt}]

% Nodes
\node[vertex, label=left:$a$]  (a) at (0,0) {};
\node[vertex, label=above:$b$] (b) at (2.5,1.5) {};
\node[vertex, label=below:$c$] (c) at (2.5,-1.5) {};
\node[vertex, label=right:$d$] (d) at (5,0) {};

\draw[->] (a) -- (b) node[midway, above left] {$\alpha_1,\beta_1$};
\draw[->] (b) -- (d) node[midway, above right] {$\alpha_2,\beta_2$};
\draw[->] (b) -- (c) node[midway, right] {$\alpha_3,\beta_3$};
\draw[->] (a) -- (c) node[midway, below left] {$\alpha_4,\beta_4$};
\draw[->] (c) -- (d) node[midway, below right] {$\alpha_5,\beta_5$};

\end{tikzpicture}
\caption{Classical Braess network configuration $N/N^+$, where $N = N^+ - (b,c)$, with arbitrary link parameters.}
\label{fig:braess-network}
\end{figure}

Let us consider the classical network configuration shown in Figure~\ref{fig:braess-network}.
This configuration is denoted by $N/N^+$, where $N = N^+ - (b,c)$. We assume that linear travel-time functions of links have arbitrary parameters.

Throughout this paper, we follow the notation used in Chapter 2 of \cite{VZverovich2021}. 
Let us denote the links in $N/N^+$ as follows:
\[
e_1=(a,b),\quad e_2=(b,d),\quad e_3=(b,c),\quad
e_4=(a,c),\quad e_5=(c,d).
\]
For each link $e_i$, 
the free-flow travel time  \(\alpha_i\) and the delay parameter \(\beta_i\) determine its travel time:
\[
T_i=\alpha_i+\beta_i f_i,
\]
where \(f_i\) denotes the flow on link \(e_i\). We assume throughout that for all $i$,
\[
\alpha_i\geq 0
\qquad\text{and}\qquad
\beta_i>0.
\]
We also adopt the following shorthand:
\[
\alpha_{ij}=\alpha_i+\alpha_j,\qquad
\bar{\alpha}=\alpha_4-\alpha_{13},\qquad
\hat{\alpha}=\alpha_2-\alpha_{35}.
\]

Moreover, let \(Q > 0\) denote the total flow in $N/N^+$, that is,
\[
Q = f_{1} + f_{4} = f_{2} + f_{5},
\]
where
\[
f_i\geq0\quad\forall i,\qquad f_{1}=f_{2}+f_{3},\qquad f_{5}=f_{3}+f_{4}.
\]
We set \(f_3=0\) in \(N\), since the link \((b,c)\) is absent.

The following definitions are well known:

\begin{definition}
A network with one origin and one destination is said to be at \emph{equilibrium} if
\begin{enumerate}
    \item[(a)] The travel time on paths with non-vanishing flow is the same (it is denoted by $T_{Eq}$) and
    \item[(b)] The travel time on paths with no flow is at least $T_{Eq}$.
\end{enumerate}
\end{definition}

\begin{definition}
Braess' paradox \emph{is said to occur} in the network configuration \(N/N^+\) for a given total flow \(Q>0\) if
\[
T_{Eq}^+ > T_{Eq},
\]
where \(T_{Eq}\) and \(T_{Eq}^+\) are travel times at equilibria in \(N\) and \(N^+\), respectively.
\end{definition}

By Theorem 2.6 in \cite{VZverovich2021}, Braess' paradox occurs in the classical network \(N/N^+\) for some \(Q>0\) if and only if the following inequalities hold:
\begin{align}
\beta_1 \beta_5 > \beta_2 \beta_4,\;\, \label{eq:BP1} \\
\beta_1 \hat{\alpha} + \beta_2 \bar{\alpha} > 0, \label{eq:BP2} \\
\beta_4 \hat{\alpha} + \beta_5 \bar{\alpha} > 0. \label{eq:BP3}
\end{align}

\section{Hierarchical networks} \label{SectionHN}

In this paper, we let \(\mathbb{N} = \{1, 2, 3, \dots\}\) represent the set of positive integers.
Let us recall the definition of a hierarchical network in \cite{Avineri2016}, which we denote by \(V/V^+\). 
In the classical network configuration, let \(
h_1, h_2, \ldots, h_5 \in \mathbb{N}
\) be the respective hierarchies of the links $e_1=(a,b)$, $e_2=(b,d)$, $e_3=(b,c)$, $e_4=(a,c)$, and $e_5=(c,d)$.
The travel time on each link \(e_i\) in \(V/V^+\) is given by
\[
T_i = \alpha k^{h_i-1} + \beta l^{h_i-1} f_i,
\]
where the parameters $k>0$ and $l>0$ are called {\em measures of hierarchy}.
The traffic flow parameters $\alpha$ and $\beta$ are
typical for roads with the highest hierarchical level ($h_i=1$).

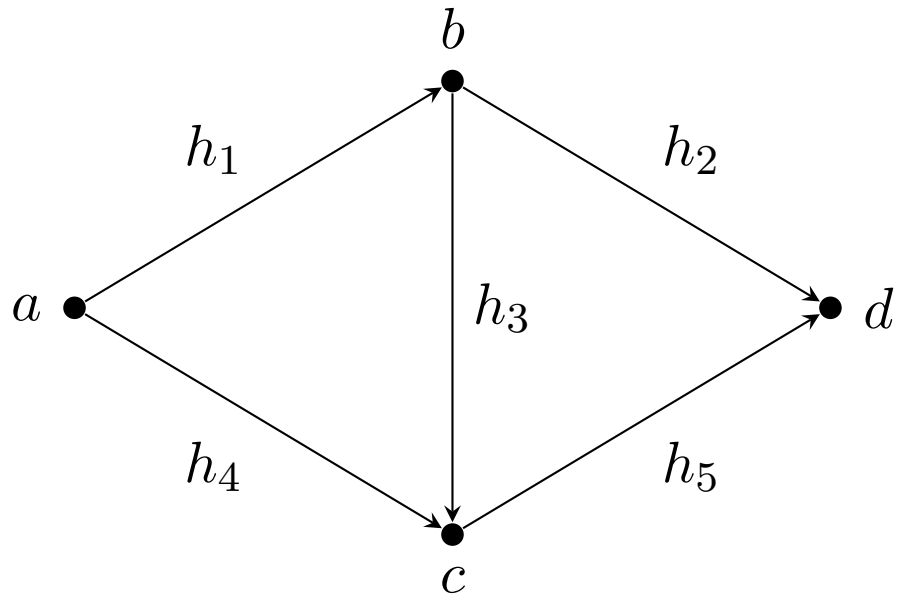
\begin{figure}[H]
\centering
\begin{tikzpicture}[
  >=stealth,
  node distance=3cm,
  vertex/.style={circle, fill=black, inner sep=1.5pt, label distance=3pt}
]

% Nodes
\node[vertex, label=left:$a$]  (a) at (0,0) {};
\node[vertex, label=above:$b$] (b) at (2.5,1.5) {};
\node[vertex, label=below:$c$] (c) at (2.5,-1.5) {};
\node[vertex, label=right:$d$] (d) at (5,0) {};

% Links
\draw[->] (a) -- (b) node[midway, above left] {$h_1$};
\draw[->] (b) -- (d) node[midway, above right] {$h_2$};
\draw[->] (b) -- (c) node[midway, right] {$h_3$};
\draw[->] (a) -- (c) node[midway, below left] {$h_4$};
\draw[->] (c) -- (d) node[midway, below right] {$h_5$};

\end{tikzpicture}
\caption{Hierarchical network configuration \(V/V^+\):
$h_i \in \mathbb{N}$,
$\alpha>0$, $\beta>0$, $k>0$, $l>0$.}
\label{fig:hierarch-braess-network}
\end{figure} 

Now, we have the following for each link \(e_i\):
\begin{align*}
\alpha_i = \alpha k^{h_i-1}
\;\;\mbox{ and }\;\;
\beta_i = \beta l^{h_i-1}.
\end{align*} 
We assume that \(\alpha>0\) and \(\beta>0\). The case \(\alpha=0\) is excluded for simplicity of presentation because it would imply that all \(\alpha_i=0\); that is, \eqref{eq:BP2} is not satisfied and Braess' paradox cannot occur. So we can compute:
\begin{align*}
\bar{\alpha}
    &= \alpha_4 - \alpha_{13} \\
    &= \alpha_4 - (\alpha_1+\alpha_3) \\
    &= \alpha k^{h_4-1} - \left(\alpha k^{h_1-1}+\alpha k^{h_3-1}\right) \\
    &= \frac{\alpha}{k}(k^{h_4}-k^{h_1} - k^{h_3})
\end{align*}
and
\begin{align*}
\hat{\alpha}
    &= \alpha_2 - \alpha_{35} \\
    &= \alpha_2 - (\alpha_3+\alpha_5) \\
    &= \alpha k^{h_2-1} - \left(\alpha k^{h_3-1}+\alpha k^{h_5-1}\right) \\
    &= \frac{\alpha}{k}(k^{h_2}-k^{h_3} - k^{h_5}).
\end{align*}

Since Theorem 2.6 gives the necessary and sufficient conditions for Braess' paradox to occur, we substitute the hierarchical link parameters into inequalities 
\eqref{eq:BP1}--\eqref{eq:BP3} to determine whether they can simultaneously hold.

For the first inequality \eqref{eq:BP1}, we obtain
\[\begin{aligned}
&\beta_1\beta_5 - \beta_2\beta_4 > 0 \\
    \Longleftrightarrow\quad &\!\left(\frac{\beta}{l}l^{h_1}\right)\left(\frac{\beta}{l}l^{h_5}\right)
       - \left(\frac{\beta}{l}l^{h_2}\right)\left(\frac{\beta}{l}l^{h_4}\right) > 0 \\
    \Longleftrightarrow\quad &\frac{\beta^2}{l^2}l^{h_1+h_5}
       - \frac{\beta^2}{l^2}l^{h_2+h_4} > 0 \\
    \Longleftrightarrow\quad &  \frac{\beta^2}{l^2}\left(l^{h_1+h_5}-l^{h_2+h_4}\right) >0 \\
    \Longleftrightarrow\quad & l^{h_1+h_5}-l^{h_2+h_4} >0 \\
    \Longleftrightarrow\quad & l^{h_1+h_5} > l^{h_2+h_4}.  
\end{aligned}\]

\noindent When \(l=1\), both sides equal $1$, so the strict inequality fails. Thus, Braess’ paradox cannot occur. Assuming $l>1$, the exponential function $l^h$ is strictly increasing in $h$, and hence this simplifies to
\[
h_1+h_5>h_2+h_4. \tag{1*}
\]
For \(0<l<1\), this condition instead becomes \(h_1+h_5<h_2+h_4\).

For the second inequality \eqref{eq:BP2}, we have
\[
\begin{aligned}
&\, \beta_1 \hat{\alpha} + \beta_2 \bar{\alpha} > 0\\
\Longleftrightarrow\quad&
\beta l^{h_1 - 1}\left(\frac{\alpha}{k}\left(k^{h_2}- k^{h_3} - k^{h_5}\right)\right)
+
\beta l^{h_2 - 1}\left(\frac{\alpha}{k}\left(k^{h_4} - k^{h_1} - k^{h_3}\right)\right)
>0\\
\Longleftrightarrow\quad&
\frac{\alpha\beta}{kl}
\left[
l^{h_1}\left(k^{h_2}-k^{h_3}-k^{h_5}\right)
+
l^{h_2}\left(k^{h_4}-k^{h_1}-k^{h_3}\right)
\right]
>0.
\end{aligned} \]
Since \(\alpha, \beta, k,l >0\), this can be simplified to
\[
l^{h_1}\left(k^{h_2}-k^{h_3}-k^{h_5}\right)
+
l^{h_2}\left(k^{h_4}-k^{h_1}-k^{h_3}\right)
>0. \tag{2*}
\]

Let us consider the third inequality \eqref{eq:BP3}:
\[
\begin{aligned}
&\beta_4 \hat{\alpha} + \beta_5 \bar{\alpha} > 0 \\
\Longleftrightarrow\quad&
\beta l^{h_4 - 1}\left(\frac{\alpha}{k}\left(k^{h_2}-k^{h_3} - k^{h_5}\right)\right)
+
\beta l^{h_5 - 1}\left(\frac{\alpha}{k}\left(k^{h_4} - k^{h_1} - k^{h_3}\right)\right)
>0\\
\Longleftrightarrow\quad&
\frac{\alpha\beta}{kl}
\left[
l^{h_4}\left(k^{h_2}-k^{h_3}-k^{h_5}\right)
+
l^{h_5}\left(k^{h_4}-k^{h_1}-k^{h_3}\right)
\right]
>0. 
\end{aligned}
\]
We can again simplify this to
\[
l^{h_4}\left(k^{h_2}-k^{h_3}-k^{h_5}\right)
+
l^{h_5}\left(k^{h_4}-k^{h_1}-k^{h_3}\right) > 0. \tag{3*}
\]

Thus, we have proved the following theorem, which provides necessary and sufficient conditions for Braess' paradox to occur in the hierarchical network configuration \(V/V^+\).

\begin{theorem} \label{thm:Braess_HierNetw}
Assume that $l>1$ and $\alpha,\beta,k > 0$. Then Braess' paradox occurs in the hierarchical network configuration \(V/V^+\) (see Figure~\ref{fig:hierarch-braess-network}) for some total flow \(Q>0\) if and only if the following inequalities hold:
\[
h_1+h_5>h_2+h_4, \tag{1*}\label{eq:HBP1}
\]
\[
l^{h_1}\left(k^{h_2}-k^{h_3}-k^{h_5}\right)
+
l^{h_2}\left(k^{h_4}-k^{h_1}-k^{h_3}\right)
>0, \tag{2*}\label{eq:HBP2}
\]
\[
l^{h_4}\left(k^{h_2}-k^{h_3}-k^{h_5}\right)
+
l^{h_5}\left(k^{h_4}-k^{h_1}-k^{h_3}\right) > 0. \tag{3*}\label{eq:HBP3}
\]
\end{theorem}

As noted above, if $0<l<1$, then 
Theorem \ref{thm:Braess_HierNetw}
holds with the first inequality reversed, namely, $h_1+h_5 < h_2+h_4$.
Also, we can immediately see that if \(k=1\), then
$$k^{h_2}-k^{h_3}-k^{h_5}=-1<0 \quad\mbox{and}\quad k^{h_4}-k^{h_1}-k^{h_3}=-1<0,$$ 
so neither \eqref{eq:HBP2} nor \eqref{eq:HBP3} can be satisfied. Therefore, Braess' paradox cannot occur in \(V/V^+\) if \(k=1\).
Moreover, the theorem has the following consequence:

\begin{corollary} \label{cor:k=l}
    Suppose \(k=l>1\). Then Braess' paradox cannot occur in the hierarchical network configuration \(V/V^+\) for any total flow \(Q>0\) and hierarchy parameters \(h_1, h_2, \ldots, h_5\). 
\end{corollary}

\begin{proof}
Suppose \(k=l>1\). It suffices to show that conditions
\eqref{eq:HBP1} and \eqref{eq:HBP2} cannot hold simultaneously.
Indeed, if \eqref{eq:HBP1} holds, then
\(h_1+h_5>h_2+h_4\), so
\(k^{h_1+h_5}>k^{h_2+h_4}\).
Setting \(l=k\), the left-hand side of \eqref{eq:HBP2} becomes
\[
\begin{aligned}
&k^{h_1}(k^{h_2}-k^{h_3}-k^{h_5})
+k^{h_2}(k^{h_4}-k^{h_1}-k^{h_3})\\
&\qquad =
k^{h_2+h_4}-k^{h_1+h_5}
-k^{h_3}(k^{h_1}+k^{h_2})<0.
\end{aligned}
\]
Thus, \eqref{eq:HBP2} fails whenever \eqref{eq:HBP1} holds.
By Theorem~\ref{thm:Braess_HierNetw}, Braess' paradox cannot occur
for any \(Q>0\).
\end{proof}

We can also derive a condition under which Braess' paradox is impossible for any choice of hierarchy parameters in a given set.

\begin{corollary} \label{cor:hi_range}
Let $h_i$, $i=1,2,\ldots,5$, be positive integers whose values lie in a set $S$ of range $r := \max{S} - \min{S}>0$.
If
\[
1<k\leq 2^{1/r}\!,
\]
then Braess' paradox cannot occur in \(V/V^+\) for any choice of
\(h_1,h_2,\ldots,h_5\) from the set $S$, any \(l>1\), or any total flow
\(Q>0\).
\end{corollary}

\begin{proof}
Let \(m\) and \(M\) be the smallest and largest values in the set $S$, so that \(r=M-m\). Since \(1<k\leq 2^{1/r}\), we have
\(\log_k 2\geq r\). Thus, for every choice of the hierarchy parameters \(h_i\) from 
$S$,
\[
h_2-\log_k 2\leq h_2-r\leq M-r=m
\leq\min\{h_3,h_5\}
\]
and, similarly,
\[
h_4-\log_k 2\leq\min\{h_1,h_3\}.
\]
Consequently,
\[
k^{h_2}\leq 2k^{\min\{h_3,h_5\}}\leq k^{h_3}+k^{h_5}
\qquad\mbox{and}\qquad
k^{h_4}\leq 2k^{\min\{h_1,h_3\}}\leq k^{h_1}+k^{h_3}.
\]
Both brackets in condition~\eqref{eq:HBP2} are therefore non-positive.
Since $l^{h_1}>0$ and 
$l^{h_2}>0$, condition~\eqref{eq:HBP2} cannot hold for
any \(l>1\). By Theorem~\ref{thm:Braess_HierNetw}, the paradox cannot
occur for any \(Q>0\).
\end{proof}

From this point onward, we assume that \(k>1\) and \(l > 1\), with the exception of Example~\(2\). The assumption reflects the idea that hierarchy 1 represents motorways, while increasing hierarchy values correspond to lower road classes (i.e.~A, B, etc.) and, consequently, longer travel times.

\section{Braess' Paradox
in \(V/V^+\)} \label{Section4}

We now investigate whether Braess' paradox can occur in the hierarchical network configuration \(V/V^+\) shown in Figure~\ref{fig:hierarch-braess-network}; that is, whether there exist values of \( k,l,h_1, h_2, \ldots, h_5\)  satisfying \eqref{eq:HBP1}--\eqref{eq:HBP3}.

First, we will develop a systematic way to find hierarchical network configurations \(V/V^+\)\!, where the paradox occurs given a fixed parameter $k>1$.

\subsection{Conditions for the occurrence of Braess' paradox in \(V/V^+\)}

We define the following:
\[
A := k^{h_2}-k^{h_3}-k^{h_5},
\]
\[
B := k^{h_4}-k^{h_1}-k^{h_3}
\]
and
\[
g_1(l) := l^{h_1}A + l^{h_2}B,
\]
\[
g_2(l) := l^{h_4}A +l^{h_5}B.
\]
Note that \(g_1(l) > 0\) is equivalent to condition~\eqref{eq:HBP2}, and \(g_2(l)>0\) is equivalent to condition~\eqref{eq:HBP3}.

The following theorem gives sufficient conditions for Braess' paradox to occur in the hierarchical network \(V/V^+\)\!.

\begin{theorem}\label{thm:Braess_existence}
Fix $k>1$ and hierarchy parameters $h_1,h_2,\ldots,h_5$. Suppose that condition~\eqref{eq:HBP1} holds and that either of the following conditions is satisfied:
\begin{enumerate}
    \item[(1)] $A>0$, $h_1>h_2$, and $h_4>h_5$;
    \item[(2)] $B>0$, $h_1<h_2$, and $h_4<h_5$.
\end{enumerate}
Then there exists an $\hat{l}>1$ such that, for every $l\geq\hat{l}$, Braess' paradox occurs in \(V/V^+\) for some total flow \(Q>0\).
\end{theorem}

\begin{proof} It suffices to prove Case~1, since the proof of Case~2 is analogous. Suppose that Case~1 holds. Since $A>0$, $h_1>h_2$, and $h_4>h_5$,
the dominant terms of $g_1$ and $g_2$ are $Al^{h_1}$ and
$Al^{h_4}$, respectively. Therefore,
\[
\lim_{l\to\infty}g_1(l)=+\infty,
\qquad
\lim_{l\to\infty}g_2(l)=+\infty.
\]
Thus, there exist $l_1>1$ and 
$l_2>1$ such that
\[
g_1(l)>0 \quad\text{for every}\quad
l\geq l_1,
\qquad
g_2(l)>0 \quad\text{for every}\quad 
l\geq l_2.
\]
Set
\[
\hat{l}:=\max\{l_1,l_2\}.
\]
Then, for every $l\geq\hat{l}$, we have $g_1(l)>0$ and $g_2(l)>0$, so conditions~\eqref{eq:HBP2} and~\eqref{eq:HBP3} hold. Since \eqref{eq:HBP1} is assumed, Theorem~\ref{thm:Braess_HierNetw} implies
that Braess' paradox occurs in \(V/V^+\) for some total flow \(Q>0\).
\end{proof}

Note that in Case 1, condition~\eqref{eq:HBP1} can always be satisfied by choosing \(h_1\) sufficiently large. Similarly, in Case 2, condition~\eqref{eq:HBP1} can be satisfied by choosing \(h_5\) sufficiently large.

Theorem~\ref{thm:Braess_existence} also yields probabilistic bounds when the hierarchy parameters are randomly distributed. 
Let us denote
\[
\begin{aligned}
    \phi_A &:= \max\{h_3,h_5\},
    &\qquad
    \phi_B &:= \max\{h_1,h_3\}, \\
    \psi_A &:= \min\{h_3,h_5\},
    &
    \psi_B &:= \min\{h_1,h_3\}.
\end{aligned}
\]

%\medskip 

\begin{corollary} \label{cor:Probability_dist_hi}
Let $D_i$ be discrete distributions supported on the positive integers, and suppose that \(h_i\sim D_i\) for \(i=1,2,\ldots,5\). No independence is assumed.
For any fixed $k>1$, let \(P_k\) denote the probability that Braess' paradox occurs in \(V/V^+\) for some \(Q>0\) and some $l>1$, 
given the parameters $h_i$.
Then 
$$
\begin{aligned}
&P_k \;\geq\; \mathbb{P}\!\left[
    \phi_A < h_2-\log_k 2
    \;\; {\rm and}\;\;
    h_1-h_2 > h_4-h_5 > 0
\right]
\\
&\hspace{4.5cm}
+\;\mathbb{P}\!\left[
    \phi_B < h_4-\log_k 2
    \;\; \text{\rm and}\;\;
    h_5-h_4 > h_2-h_1 > 0
\right] 
\end{aligned}
\vspace*{-0.1cm}
$$
and
\vspace*{-0.1cm}
$$
P_k \leq
1 - \mathbb{P}\!\left[
        h_1+h_5 \leq h_2+h_4 
        \;\;\text{\rm or}\;\;
        (h_2 - \log_k 2 \leq \psi_A \;\; \text{\rm and}\;\; h_4 - \log_k 2 \leq \psi_B)
\right].
$$
\end{corollary}

The proof is given in Appendix~\ref{app:Probability_dist_hi_proof}.
Although these bounds can yield a wide range of possible values for $P_k$, the lower bounds provide a useful guarantee on the minimum probability of Braess’ paradox occurring under the specified distribution. In some cases, the upper bounds are also informative. For example, the upper bound for $P_{1.1}$ is 0.01129, assuming that the hierarchy parameters are independently and uniformly distributed on \(\{1,2,\ldots,10\}\). This bound was evaluated by enumerating all \(10^5\) equally likely configurations of the hierarchy parameters.

\subsection{Examples of Braess' paradox in \(V/V^+\)}

We now present an explicit example of Braess' paradox in a hierarchical network. Since the link $(b,c)$ is the newly added road, we begin by restricting our attention to the case $h_3=1$, so that it belongs to the highest/fastest class of roads. This simplifies the search for suitable parameters. However, as shown in Theorem~\ref{thm:Braess_existence}, for any $k>1$, Braess' paradox may occur for arbitrary values of $h_3$, provided that condition~\eqref{eq:HBP1} holds, $A>0$, $h_1>h_2$, $h_4>h_5$ and $l$ is sufficiently large.

\medskip

\begin{example} Consider the hierarchical network \(V_1/V_1^+\)\! with
\begin{align*}
&\alpha>0,\; \beta>0,\;
k=3, \; l=9, \\ &h_1=5, \; h_2=2, \; h_3=1, \; h_4=3, \; h_5=1.
\end{align*}
\end{example}

Let us verify that conditions~\eqref{eq:HBP1}--\eqref{eq:HBP3} are satisfied in $V_1/V_1^+$.
Clearly,
$h_1+h_5=6>5=h_2+h_4,$
so \eqref{eq:HBP1} holds.
Further,
\begin{align*}
&l^{h_1}\left(k^{h_2}-k^{h_3}-k^{h_5}\right)
+
l^{h_2}\left(k^{h_4}-k^{h_1}-k^{h_3}\right) = 159,408 > 0,
\end{align*}
so \eqref{eq:HBP2} holds.
Finally,
\begin{align*}
&l^{h_4}\left(k^{h_2}-k^{h_3}-k^{h_5}\right)
+
l^{h_5}\left(k^{h_4}-k^{h_1}-k^{h_3}\right)=216>0,
\end{align*}
so \eqref{eq:HBP3} holds.

Therefore, by Theorem~\ref{thm:Braess_HierNetw}, for any
$\alpha>0$ and $\beta>0$ there exists a total flow \(Q>0\) for which Braess' paradox occurs in $V_1/V_1^+$.
Since conditions~\eqref{eq:HBP2} and~\eqref{eq:HBP3} remain satisfied for all \(l \ge 9\), Braess' paradox also occurs for every \(l \ge 9\).

In the previous example, we have $k>1$, $l>1$, and $h_3=1$. 
We now demonstrate another example in which $0<k<1$, $0<l<1$, and $h_3\ge 3$.

\medskip

\begin{example} Consider the hierarchical network \(V_2/V_2^+\)\! with
\begin{align*}
%&\alpha=27,\; \beta=1,\;
&\alpha>0,\; \beta>0,\;
k={1\over 3}, \; l={9\over 10}, \\ &h_1=2, \; h_2=4, \; h_3\ge 3, \; h_4=1, \; h_5=2.
\end{align*}
\end{example}

\begin{proposition}
Braess’ paradox occurs in the hierarchical network $V_2$ for some total flow \(Q>0\) when adding the link $e_3$ of any hierarchy $h_3\ge 3$.
\end{proposition}

The proof is left to the reader. Note that since $l<1$, we must reverse the inequality \eqref{eq:HBP1} in Theorem~\ref{thm:Braess_HierNetw}.

\section{Monte Carlo simulations} \label{Section5}

We use Monte Carlo simulation to estimate the probability that the inequalities in Theorem~\ref{thm:Braess_HierNetw} are simultaneously satisfied under different distributions of the network parameters. These probabilities describe the frequency with which a sampled network admits the paradox for some total flow \(Q>0\) under the stated parameter distributions. They should not be interpreted as estimates of its prevalence in real road networks.

These probabilities differ from \(P_k\) in Corollary~\ref{cor:Probability_dist_hi}, which fixes \(k>1\) and measures the probability, over the random hierarchy parameters, that there exist \(l>1\) and \(Q>0\) for which the paradox occurs. In the simulations below, \(k\), \(l\) and the hierarchy parameters are sampled, and we test whether the resulting network admits the paradox for some \(Q>0\).

All simulations are implemented in Python 3.12.8 using NumPy version 1.26.4, with Jupyter Notebook, through the Anaconda distribution.
Random samples are generated using NumPy's PCG64 random-number generator, initialised with seed 101 via  \texttt{numpy.random.default\_rng}.

The probabilities reported below are Monte Carlo estimates and are therefore subject to sampling uncertainty. Small differences between estimates should therefore be interpreted with care. For each sample, the parameters are drawn from their respective distributions, and the inequalities in Theorem~\ref{thm:Braess_HierNetw} are evaluated to determine whether Braess' paradox occurs. The estimated probability is calculated as the proportion of samples for which all the required inequalities are satisfied. The estimated Monte Carlo standard error of an estimated probability
$\widehat{p}$ is
\[
\sqrt{\frac{\widehat{p}(1-\widehat{p})}{N}}.
\]

For the rest of this paper, every simulation will use \(N=10^7\) independent samples. Tabulated probability estimates will be rounded to two significant figures, while estimates accompanied by confidence intervals will be rounded to the same decimal place as their margins of error. Margins of error are rounded to two significant figures. Moreover, any estimated probability below $10^{-5}$, including those for which no occurrences are observed, will be denoted as \(<10^{-5}\)\!. However, this notation refers to the estimated probability and does not establish an upper bound on the true probability. The main probability estimates discussed will be accompanied by approximate $95\%$ confidence intervals, calculated as the estimate $\pm 1.96$ times its estimated Monte Carlo standard error.

The hierarchy parameters \(h_1, h_2, \ldots, h_5\) are assumed to be independent and identically distributed according to the specified discrete distribution. 
Unless stated otherwise, the hierarchy measures $k$ and $l$ 
are mutually independent and independent of the hierarchy parameters, and are identically distributed according to the specified continuous distribution.

To facilitate comparison between distributions, their parameters are chosen so that they have the same mean. For $k$ and $l$, the uniform distribution is supported on $[1,10]$, while the exponential, Erlang, log-normal and beta distributions are shifted by 1 and parametrised to have mean $5.5$ and support contained in $[1,\infty)$. Specifically, the exponential distribution has scale $4.5$, and the Erlang distributions with shape parameters $2$, $4$, and $16$ have scale parameters $2.25$, $1.125$, and $0.28125$, respectively. The log-normal distribution has parameters $\mu=0$ and $\sigma=1.7344$ and is also shifted by 1. For the beta distribution, let $X_k$ and $X_l$ be independent random variables, each following a \(\operatorname{Beta}(6,2)\) distribution, and define \(k=1+6X_k\) and \(l=1+6X_l\).

The hierarchy parameters are integer-valued and are drawn from one of five distributions, each with mean $5.5$:

\vspace*{-0.2cm}
\begin{enumerate}
    \item $\mathrm{Uniform}(1,10)$;
    \vspace*{-0.2cm}\item $1 + \mathrm{Binomial}(9,0.5)$;
    \vspace*{-0.2cm}\item $1 + \mathrm{NegativeBinomial}(9,2/3)$;
    \vspace*{-0.2cm}\item $1 + \mathrm{Poisson}(4.5)$;
    \vspace*{-0.2cm}\item $\mathrm{Geometric}(1/5.5)$.
\end{enumerate}

\vspace*{-0.2cm}
The geometric distribution is assumed to have support $\{1,2,\ldots\}$, whereas the negative binomial distribution counts the number of failures before a specified number of successes.

\vspace*{0.3cm}
\begin{table}[H]
\centering
\renewcommand{\arraystretch}{1.6}
\resizebox{\textwidth}{!}{%
\small
\begin{tabular}{|c||c|c|c|c|c|}
\hline
\raisebox{-3mm}{
\textbf{Distributions}} 
& \multicolumn{5}{c|}{\textbf{Distribution of Hierarchy Parameters, $h_i$:}} \\
\cline{2-6}
\textbf{of $k$ and $l$:} & \textbf{Uniform} & \textbf{Binomial} & \textbf{Neg-Binomial} & \textbf{Poisson} & \textbf{Geometric}\\
\hline\hline
\textbf{Uniform(1,10)} & 0.0051 & 0.0019 & 0.0043 & 0.0034 & 0.0051 \\
\hline
\textbf{Exponential} & 0.0072 & 0.0028 & 0.0060 & 0.0049 & 0.0069 \\
\hline
\textbf{Erlang-2} & 0.0047 & 0.0016 & 0.0041 & 0.0032 & 0.0050 \\
\hline
\textbf{Erlang-4} & 0.0025 & 0.00065 & 0.0023 & 0.0016 & 0.0031 \\
\hline
\textbf{Erlang-16} & 0.00043 & 0.000037 & 0.00038 & 0.00021 & 0.00084 \\
\hline
\textbf{Log-normal(0,1.7344)} & 0.0085 & 0.0035 & 0.0066 & 0.0055 & 0.0074 \\
\hline
\textbf{Beta(6,2)} & 0.00032 & 0.000034 & 0.00028 & 0.00016 & 0.00063 \\
\hline
\end{tabular}%
}
\caption{Estimated probability of Braess' paradox in hierarchical networks under different distributions of the network parameters.}
\label{tab:network-parameters}
\end{table}

From Table~\ref{tab:network-parameters}, we can see that the probability of Braess' paradox occurring in hierarchical networks \(V/V^+\) is relatively small across all of the distributions considered, with estimated probabilities below \(1\%\) in every case. However, the choice of distribution for $k$ and $l$ has a noticeable effect. The largest probabilities occur for the log-normal and exponential distributions, while the Erlang-16 and Beta(6,2) distributions produce substantially smaller probabilities. In particular, as the Erlang shape parameter increases, the estimated probability decreases, suggesting that greater concentration of $k$ and $l$ around their common mean makes the paradox less likely.

\subsection{Mean}

Let us now study the effect of shifts in the means of the network parameter distributions on the probability of Braess' paradox occurring in hierarchical networks.

In the following table, we again assume that the hierarchy parameters are independent and uniformly distributed on \(\{1, 2, \ldots, 10\}\). 
Let us investigate the effect of changes in the means of $k$ and $l$. Both parameters are independently uniformly distributed, with their supports chosen to give the required means. In particular, for a desired mean value $\mu_k$ of $k$, we take $k$ to be uniformly distributed on \([\mu_{k} - 0.5, \mu_{k} + 0.5]\), and similarly for a mean value $\mu_l$ of $l$,  we take $l$ to be uniformly distributed on \([\mu_{l} - 0.5, \mu_{l} + 0.5]\). We consider values of $\mu_k$ and $\mu_l$ between $1.5$ and $9.5$. 

\medskip

\begin{table}[H]
\centering
\renewcommand{\arraystretch}{1.3}
\resizebox{\textwidth}{!}{%
\begin{tabular}{|c||c|c|c|c|c|c|c|c|c|}
\hline
\raisebox{-1mm}{\textbf{Mean of the}}
& \multicolumn{9}{c|}{\textbf{Mean of the distribution of \textit{l}}} \\
\cline{2-10}
\textbf{distribution of \textit{k}}& \textbf{1.5} & \textbf{2.5} & \textbf{3.5} & \textbf{4.5}
& \textbf{5.5} & \textbf{6.5} & \textbf{7.5} & \textbf{8.5} & \textbf{9.5} \\
\hline\hline

\textbf{1.5}
& 0.0019 & 0.00079 & 0.0024 & 0.0036 & 0.0043 & 0.0048 & 0.0052 & 0.0055 & 0.0057 \\
\hline

\textbf{2.5}
& 0.012 & 0.00023 & 0.0010 & 0.0029 & 0.0045 & 0.0057 & 0.0065 & 0.0074 & 0.0080 \\
\hline

\textbf{3.5}
& 0.019 & 0.0026 & $<10^{-5}$ & 0.00028 & 0.0013 & 0.0025 & 0.0035 & 0.0043 & 0.0047 \\
\hline

\textbf{4.5}
& 0.023 & 0.0065 & 0.00060 & $<10^{-5}$ & 0.000039 & 0.00047 & 0.0013 & 0.0020 & 0.0028 \\
\hline

\textbf{5.5}
& 0.025 & 0.0099 & 0.0025 & 0.000094 & $<10^{-5}$ & $<10^{-5}$ & 0.00013 & 0.00057 & 0.0012 \\
\hline

\textbf{6.5}
& 0.026 & 0.012 & 0.0047 & 0.00089 & $<10^{-5}$ & $<10^{-5}$ & $<10^{-5}$ & 0.000020 & 0.00020 \\
\hline

\textbf{7.5}
& 0.027 & 0.014 & 0.0065 & 0.0023 & 0.00026 & $<10^{-5}$ & $<10^{-5}$ & $<10^{-5}$ & $<10^{-5}$ \\
\hline

\textbf{8.5}
& 0.028 & 0.016 & 0.0079 & 0.0036 & 0.0010 & 0.000051 & $<10^{-5}$ & $<10^{-5}$ & $<10^{-5}$ \\
\hline

\textbf{9.5}
& 0.029 & 0.017 & 0.0089 & 0.0049 & 0.0021 & 0.00041 & $<10^{-5}$ & $<10^{-5}$ & $<10^{-5}$ \\
\hline

\end{tabular}%
}
\caption{Estimated probability of Braess' paradox for different means of $k$ and $l$.}
\label{tab:braess_kl_mean}
\end{table}

Table~\ref{tab:braess_kl_mean} suggests that the estimated probability of Braess' paradox is strongly associated with the relative means of $k$ and $l$. In particular, the highest probabilities occur when the mean of $k$ is large relative to the mean of $l$. For example, when the mean of $l$ is fixed at $1.5$, increasing the mean of $k$ from $1.5$ to $9.5$ increases the estimated probability from $0.001911 \pm 2.7 \times 10^{-5}$ to $0.02877 \pm 1.0 \times 10^{-4}$. 

Conversely, when the mean of $k$ is large, increasing the mean of $l$ substantially reduces the estimated probability. For example, when the mean of $k$ is $9.5$, increasing the mean of $l$ from $1.5$ to $9.5$ reduces the estimated probability from $0.02877$ to below $10^{-5}$, with no occurrences observed in the latter case. The relationship is therefore not simply that larger values of either parameter increase the probability. Rather, the results suggest that Braess' paradox could be more likely when the distribution of $k$ has a considerably larger mean than the distribution of $l$. However, Theorem~\ref{thm:Braess_existence} shows that, for fixed $k$ and hierarchy parameters satisfying its conditions, Braess' paradox occurs for all sufficiently large $l$. This does not establish that increasing the mean of $l$ increases the overall probability when averaging over randomly sampled hierarchy parameters.

Corollary~\ref{cor:Probability_dist_hi} implies that, with the distributions of the hierarchy parameters held fixed, the probability of Braess' paradox occurring for some $l>1$ and some \(Q>0\) tends to zero as $k \to 1^+$. Indeed, its upper bound is
\[
1-\mathbb{P}\!\left[
    h_1+h_5 \leq h_2+h_4
    \ \text{or}\
    \bigl(
        h_2-\log_k 2 \leq \psi_A
        \ \text{and}\
        h_4-\log_k 2 \leq \psi_B
    \bigr)
\right].
\]
As $k\rightarrow 1^+$, we have $\log_k 2 \rightarrow \infty$. For any fixed choice of the hierarchy parameters, both inequalities in the second condition therefore eventually hold. Hence, the indicator of the event inside the probability converges pointwise to $1$. By the dominated convergence theorem, the probability tends to $1$. Consequently, the upper bound, and hence the probability of Braess' paradox occurring for some $l>1$, tends to $0$.

We now show that an admissible integer shift in the means of the hierarchy parameter distributions leaves the probability of Braess’ paradox unchanged:

\begin{lemma}[Invariance under a common shift] \label{lem:h_i-mean}
Fix $k,l>1$. For any $c\in\mathbb{Z}$ such that
$h_i+c\in\mathbb{N}$ for every $i$, Braess' paradox occurs in
\(V/V^+\) with hierarchy parameters \(h_1, h_2, \ldots, h_5\) for some total flow \(Q>0\) if and only if
it occurs in $\widetilde{V}/\widetilde{V}^+$ with hierarchy parameters
\[
    \widetilde{h}_i=h_i+c,
    \qquad i=1, 2,\ldots,5,
\]
for some total flow \(\widetilde{Q}>0\), where
\(\widetilde{Q}\) need not be equal to \(Q\).
\end{lemma}

\begin{proof}
We verify that conditions~\eqref{eq:HBP1}--\eqref{eq:HBP3} in Theorem~\ref{thm:Braess_HierNetw}
hold for \(V/V^+\) if and only if they hold for
$\widetilde{V}/\widetilde{V}^+$\!.
Let us fix $c\in\mathbb{Z}$ such that $h_i+c\in\mathbb{N}$ for every $i$.

For \eqref{eq:HBP1}, we obtain
\begin{align*}
h_1+h_5>h_2+h_4
&\iff h_1+h_5+2c>h_2+h_4+2c \\
&\iff \widetilde{h}_1+\widetilde{h}_5
    >\widetilde{h}_2+\widetilde{h}_4.
\end{align*}

For \eqref{eq:HBP2}, since $(kl)^c>0$, we have
\begin{align*}
&l^{\widetilde{h}_1}
 \left(k^{\widetilde{h}_2}
      -k^{\widetilde{h}_3}
      -k^{\widetilde{h}_5}\right)
+l^{\widetilde{h}_2}
 \left(k^{\widetilde{h}_4}
      -k^{\widetilde{h}_1}
      -k^{\widetilde{h}_3}\right)>0 \\
\iff\quad&
l^{h_1+c}
 \left(k^{h_2+c}-k^{h_3+c}-k^{h_5+c}\right)
+l^{h_2+c}
 \left(k^{h_4+c}-k^{h_1+c}-k^{h_3+c}\right)>0 \\
\iff\quad&
(kl)^c
\left[
l^{h_1}\left(k^{h_2}-k^{h_3}-k^{h_5}\right)
+l^{h_2}\left(k^{h_4}-k^{h_1}-k^{h_3}\right)
\right]>0 \\
\iff\quad&
l^{h_1}\left(k^{h_2}-k^{h_3}-k^{h_5}\right)
+l^{h_2}\left(k^{h_4}-k^{h_1}-k^{h_3}\right)>0.
\end{align*}

Similarly, for \eqref{eq:HBP3}, we obtain
\begin{align*}
&l^{\widetilde{h}_4}
 \left(k^{\widetilde{h}_2}
      -k^{\widetilde{h}_3}
      -k^{\widetilde{h}_5}\right)
+l^{\widetilde{h}_5}
 \left(k^{\widetilde{h}_4}
      -k^{\widetilde{h}_1}
      -k^{\widetilde{h}_3}\right)>0 \\
\iff\quad&
l^{h_4+c}
 \left(k^{h_2+c}-k^{h_3+c}-k^{h_5+c}\right)
+l^{h_5+c}
 \left(k^{h_4+c}-k^{h_1+c}-k^{h_3+c}\right)>0 \\
\iff\quad&
(kl)^c
\left[
l^{h_4}\left(k^{h_2}-k^{h_3}-k^{h_5}\right)
+l^{h_5}\left(k^{h_4}-k^{h_1}-k^{h_3}\right)
\right]>0 \\
\iff\quad&
l^{h_4}\left(k^{h_2}-k^{h_3}-k^{h_5}\right)
+l^{h_5}\left(k^{h_4}-k^{h_1}-k^{h_3}\right)>0.
\end{align*}
\end{proof}

\subsection{Skewness}

We now examine how skewness in the distributions of the network parameters affects the probability of Braess' paradox in hierarchical networks.

In the following table, we again assume that the hierarchy parameters are independent and uniformly distributed on \(\{1, 2, \ldots,10\}\). Moreover, let \(X_k\) and \(X_l\) be independent random variables, each following a 
\(\mathrm{Beta}(a,b)\) distribution. We set \(k=1+mX_k\) and \(l=1+mX_l\), where \(m\) is chosen for each pair \((a,b)\) so that both $k$ and $l$ have mean $5.5$. Since \(\mathbb{E} [X_k]=\mathbb{E} [X_l]=a/(a+b)\), we set \(m=4.5(a+b)/a\). It is worth noting that this scaling of the beta distribution does not change its skewness. 
Thus, $k$ and $l$ take values in \([1,1+m]\). 

\medskip
\begin{table}[H]
\centering
\renewcommand{\arraystretch}{1.1}
\small
\begin{tabular}{|c||c|c|c|c|c|c|}
\hline
\raisebox{-3mm}{
\textbf{ Parameter \textit{a}}} & \multicolumn{6}{c|}{\textbf{Parameter \textit{b}}} \\
\cline{2-7}
 & \textbf{1} & \textbf{2} & \textbf{3} & \textbf{4} & \textbf{5} & \textbf{6} \\
\hline\hline
\textbf{1} & 0.0051 & 0.0060 & 0.0064 & 0.0066 & 0.0067 & 0.0068 \\
\hline
\textbf{2} & 0.0021 & 0.0029 & 0.0033 & 0.0036 & 0.0038 & 0.0039 \\
\hline
\textbf{3} & 0.00098 & 0.0015 & 0.0019 & 0.0021 & 0.0023 & 0.0025 \\
\hline
\textbf{4} & 0.00050 & 0.00090 & 0.0012 & 0.0013 & 0.0015 & 0.0016 \\
\hline
\textbf{5} & 0.00026 & 0.00052 & 0.00073 & 0.00088 & 0.0010 & 0.0011 \\
\hline
\textbf{6} & 0.00016 & 0.00032 & 0.00047 & 0.00059 & 0.00070 & 0.00077 \\
\hline
\end{tabular}
\caption{Estimated probability of Braess' paradox for scaled and shifted Beta($a,b$) distributions of $k$ and $l$.}
\label{tab:braess_kl_skew}
\end{table}

Table~\ref{tab:braess_kl_skew} shows a clear association between the skewness of the distributions of $k$ and $l$ and the estimated probability of Braess' paradox. For a fixed value of $a$, increasing $b$ (moving the skewness in the positive direction) generally increases the estimated probability, while for a fixed value of $b$, increasing $a$ (moving the skewness in the negative direction) decreases it. Thus, greater positive skewness in the distributions of $k$ and $l$ is associated with a higher estimated probability, whereas greater negative skewness is associated with a lower estimated probability.

This difference is particularly pronounced when the parameters are reversed. For example, the estimated probability is $0.006821 \pm 5.1 \times 10^{-5}$ for Beta(1,6), compared with only $0.0001572 \pm 7.8 \times 10^{-6}$ for Beta(6,1). Although the probability remains small in absolute terms, these results suggest that greater positive skewness for the distribution of $k$ and $l$ is associated with a higher probability of Braess' paradox. However, other properties of the distributions also change in these experiments, so the observed differences in probability cannot be attributed to skewness alone. This can be seen in the $\operatorname{Beta}(1,1)$ and $\operatorname{Beta}(6,6)$ distributions, which both have zero skewness but very different estimated probabilities of $0.005072 \pm 4.4 \times 10^{-5}$ and $0.000774 \pm 1.7 \times 10^{-5}$ respectively.

We next investigate the relationship between the skewness of the distribution of the hierarchy parameters \(h_1, h_2, \ldots, h_5\) and the probability of Braess' paradox occurring. Let us assume that $k$ and $l$ are independent and uniformly distributed on $[1,10]$. For the hierarchy parameters, let \(X_1, X_2, \ldots, X_5\) be independent random variables with $X_i\sim\operatorname{Binomial}(n,p)$, and define $h_i=X_i+1$. We consider \(p\in\{0.1,0.2,0.4,0.5,0.6,0.8,0.9\}\) and choose $n$ for each value of $p$ so that the variance of each hierarchy parameter is as close as possible to $2.25$.  Since 
\(\operatorname{Var} [h_i]=
\operatorname{Var} [X_i]=np(1-p)\), we set 
\[
n=\operatorname{round}\left(\frac{2.25}{p(1-p)}\right).
\]

Although the resulting distributions have different means, Lemma~\ref{lem:h_i-mean} de\-monstrates that a common shift of hierarchy parameters does not affect whether Braess' paradox occurs. Hence, we compare distributions with approximately matched variances, ranging from $2.16$ to $2.25$.

\medskip
\begin{table}[H]
\centering
\renewcommand{\arraystretch}{1.2}
\resizebox{\textwidth}{!}{%
\small
\begin{tabular}{|c||c|c|c|c|c|c|c|}
\hline
\textbf{Parameter $p$}
& \textbf{0.1} & \textbf{0.2} & \textbf{0.4}
& \textbf{0.5} & \textbf{0.6} & \textbf{0.8}
& \textbf{0.9} \\
\hline\hline
\textbf{Parameter $n$}
& 25 & 14 & 9 & 9 & 9 & 14 & 25 \\
\hline
\textbf{Variance, $np(1-p)$}
& 2.25 & 2.24 & 2.16 & 2.25 & 2.16 & 2.24 & 2.25 \\
\hline
\textbf{Estimated probability}
& 0.0016 & 0.0017 & 0.0017 & 0.0019
& 0.0018 & 0.0019 & 0.0019 \\
\hline
\end{tabular}%
}
\caption{Estimated probability of Braess' paradox for shifted binomial distributions of the hierarchy parameters.}
\label{tab:braess_hi_skew}
\end{table}

Table~\ref{tab:braess_hi_skew} suggests a relationship between the skewness of the binomial distribution of the hierarchy parameters and the probability of Braess' paradox. Increasing $p$ from $0.1$ to $0.5$, and thus moving the skewness from positive towards zero, slightly increases the estimated probability from $0.001561 \pm 2.4 \times 10^{-5}$ to $0.001886 \pm 2.7 \times 10^{-5}$. This suggests that right-skewed discrete distributions of \(h_1,h_2,\ldots,h_5\) may be associated with a lower estimated probability of Braess' paradox. For $p>0.5$, however, greater left-skewness does not appear to correspond to a consistently higher estimated probability.

Overall, the differences are small in absolute terms: the estimated probability for $p=0.1$ is only around $0.0003$ ($0.03$ percentage points) lower than those for $p=0.5$ and $p=0.9$. Moreover, since $n$ must be integer-valued, the variance cannot be kept exactly constant for all values of $p$. These results should therefore be interpreted as evidence of some association between right-skewness of the hierarchy parameters and a lower probability of Braess' paradox.

In conclusion, the observed relationship between skewness in the distribution of the hierarchy-scaling parameters and the probability of Braess' paradox occurring is consistent with the findings presented at the end of Chapter 2 of \cite[p.~94]{VZverovich2021}.  In both the beta and binomial cases, changing the distributional parameters also affects other properties, such as the variance and support, as discussed in the next subsection. These simulations should therefore be interpreted as evidence of an association with skewness rather than establishing skewness as the sole cause of the observed changes.

\subsection{Standard Deviation}

Let us study how the standard deviation in the distributions of the network parameters affects the probability of Braess' paradox occurring in hierarchical networks.

In the following table, we also assume that the hierarchy parameters are independent and uniformly distributed on \(\{1,2,\ldots,10\}\). We take $k$ and $l$ to be independent and uniformly distributed on \([5.5-0.5r_k,\,5.5+0.5r_k]\) and \([5.5-0.5r_l,\,5.5+0.5r_l]\), respectively, where \(r_k\) and \(r_l\) denote their respective ranges. Both 
$k$ and $l$ have expected value $5.5$.  Since the standard deviation of a continuous uniform distribution with range $r$ is \(r/\sqrt{12}\), we report the ranges rather than the standard deviations.

\medskip
\begin{table}[H]
\centering
\renewcommand{\arraystretch}{1.1}
\small
\begin{tabular}{|c||c|c|c|c|c|c|}
\hline
\raisebox{-3mm}{
\textbf{Range of $k$}} & \multicolumn{6}{c|}{\textbf{Range of $l$}} \\
\cline{2-7}
 & \textbf{1} & \textbf{2} & \textbf{3} & \textbf{5} & \textbf{7} & \textbf{9} \\
\hline\hline
\textbf{1} & $<10^{-5}$ & $<10^{-5}$ & 0.000033 & 0.00054 & 0.0019 & 0.0044 \\
\hline
\textbf{2} & $<10^{-5}$ & 0.000013 & 0.000078 & 0.00061 & 0.0019 & 0.0044 \\
\hline
\textbf{3} & 0.000017 & 0.000063 & 0.00016 & 0.00071 & 0.0020 & 0.0045 \\
\hline
\textbf{5} & 0.00033 & 0.00040 & 0.00054 & 0.0011 & 0.0023 & 0.0047 \\
\hline
\textbf{7} & 0.0010 & 0.0011 & 0.0012 & 0.0017 & 0.0028 & 0.0050 \\
\hline
\textbf{9} & 0.0015 & 0.0016 & 0.0017 & 0.0022 & 0.0031 & 0.0051 \\
\hline
\end{tabular}
\caption{Estimated probability of Braess' paradox for different ranges in the distributions of $k$ and $l$.}
\label{tab:braess_kl_sd}
\end{table}

Table~\ref{tab:braess_kl_sd} shows that increasing the range of either $k$ or $l$ is associated with an increase in the estimated probability of Braess' paradox occurring. For example, with the range of $l$ fixed at $1$, increasing the range of $k$ from $1$ to $9$ increases the estimated probability from less than $10^{-5}$ to $0.001506 \pm 2.4 \times 10^{-5}$. Similarly, with the range of $k$ fixed at $1$, increasing the range of $l$ from $1$ to $9$ increases the estimated probability from less than $10^{-5}$ to $0.004360 \pm  4.1 \times 10^{-5}$.

Over the ranges considered, variation in $l$ appears to have a greater effect than variation in $k$. One possible explanation follows from Theorem~\ref{thm:Braess_existence}: widening the range of $l$ allows larger values of $l$ to occur, which can make the conditions for Braess' paradox more likely to be satisfied. 

In the following table, we again assume that $k$ and $l$ are independent and uniformly distributed on \([1,10]\). We vary the standard deviation of the hierarchy parameters by changing the range of their distributions. In particular, for a given range $r$, the hierarchy parameters are independent and uniformly distributed on \(\{1,2,\ldots,r+1\}.\) By Lemma~\ref{lem:h_i-mean}, a common shift in the mean should not affect the results. We again use the range as the primary parameter and also report the corresponding standard deviation, which increases with the range.

\medskip
\begin{table}[H]
\centering
\renewcommand{\arraystretch}{1.3}
\resizebox{\textwidth}{!}{%
\begin{tabular}{|c||c|c|c|c|c|c|c|c|}
\hline
\textbf{Range of the Distribution of $h_i$} & \textbf{1} & \textbf{3} & \textbf{6} & \textbf{9} & \textbf{12}
& \textbf{15} & \textbf{19} & \textbf{29} \\
\hline\hline
\textbf{Standard Deviation} 
& 0.5 & 1.12 & 2.0 & 2.87 & 3.74 & 4.61 & 5.77 & 8.66 \\
\hline
\textbf{Estimated probability}
& $<10^{-5}$ & 0.00048 & 0.0028 & 0.0050 & 0.0067 & 0.0079
& 0.0090 & 0.011 \\
\hline
\end{tabular}%
}
\caption{Estimated probability of Braess' paradox for different ranges of the hierarchy parameters $h_i$.}
\label{tab:braess_hi_sd}
\end{table}

Table~\ref{tab:braess_hi_sd} suggests that increasing the range, and hence the standard deviation, of the hierarchy parameters is associated with a higher estimated probability of Braess' paradox. Increasing the range from $1$ to $29$ raises the estimate from less than $10^{-5}$ to $0.010801 \pm 6.4 \times 10^{-5}$, with no occurrences observed for the former. 

One possible explanation comes from the lower bound in Corollary~\ref{cor:Probability_dist_hi}. The two events defining this bound require particular orderings of the hierarchy parameters, together with sufficiently large differences between some of them. A wider range allows larger differences, which may make these events more likely. However, the bound concerns the existence of some $l>1$, whereas the estimates use sampled values of $l$. It therefore suggests a possible reason for the observed increase but does not establish that the probability must increase with the range.

Overall, the observed relationship between the standard deviation and the probability of Braess' paradox occurring in both sets of simulations is again consistent with the findings presented at the end of Chapter 2 of \cite[p.~94]{VZverovich2021}. The simulations indicate a positive association between the ranges of the network parameter distributions (and hence their standard deviations) and the estimated probability of Braess' paradox for the distributions and ranges considered.

\subsection{Comparison of Parameter Distributions}

We now compare how different combinations of distributions for $k$ and $l$ affect the estimated probability of Braess' paradox while holding the distribution of the hierarchy parameters fixed. We consider three fixed distributions for the hierarchy parameters: $\mathrm{Uniform}(1,10)$, $1+\mathrm{Binomial}(9,0.5)$ and $\mathrm{Geometric}(1/5.5)$. For each case, we compare the seven distributions of $k$ and $l$ listed in Table~\ref{tab:network-parameters}, using the same parameter values and shifts. This allows us to identify which combinations produce the highest estimated probabilities among those considered.

\begin{table}[H]
\centering
\renewcommand{\arraystretch}{1.6}
\resizebox{\textwidth}{!}{%
\begin{tabular}{|c||c|c|c|c|c|c|c|}
\hline
\raisebox{-3mm}{
\textbf{Distribution}} 
& \multicolumn{7}{c|}{\textbf{Distribution of $l$:}} \\
\cline{2-8}
\textbf{of $k$:} & \textbf{Uniform} & \textbf{Exponential} & \textbf{Erlang-2} & \textbf{Erlang-4} & \textbf{Erlang-16} & \textbf{Log-normal} & \textbf{Beta} \\
\hline\hline
\textbf{Uniform} & 0.0051 & 0.0074 & 0.0047 & 0.0030 & 0.0018 & 0.014 & 0.0018 \\
\hline
\textbf{Exponential} & 0.0053 & 0.0072 & 0.0048 & 0.0036 & 0.0028 & 0.013 & 0.0028 \\
\hline
\textbf{Erlang-2} & 0.0052 & 0.0076 & 0.0047 & 0.0031 & 0.0021 & 0.014 & 0.0021 \\
\hline
\textbf{Erlang-4} & 0.0049 & 0.0076 & 0.0043 & 0.0025 & 0.0013 & 0.015 & 0.0012 \\
\hline
\textbf{Erlang-16} & 0.0046 & 0.0075 & 0.0039 & 0.0019 & 0.00044 & 0.016 & 0.00037 \\
\hline
\textbf{Log-normal} & 0.0046 & 0.0055 & 0.0042 & 0.0037 & 0.0035 & 0.0085 & 0.0034 \\
\hline
\textbf{Beta} & 0.0045 & 0.0076 & 0.0039 & 0.0018 & 0.00036 & 0.016 & 0.00032 \\
\hline
\end{tabular}%
}
\caption{Estimated probability of Braess' paradox under different distributions of $k$ and $l$, with the hierarchy parameters independently and uniformly distributed on \(\{1,2,\ldots,10\}\).}
\label{tab:kl-uniform_hi}
\end{table}

First, let us assume that the hierarchy parameters are independent and uniformly distributed on \(\{1,2,\ldots,10\}\).
Table~\ref{tab:kl-uniform_hi} shows that the estimated probability of Braess' paradox is highest when $l$ follows a log-normal distribution. The largest estimated probability is $0.016026 \pm 7.8 \times 10^{-5}$, occurring when $k$ follows a beta distribution. Conversely, the lowest estimate is $0.000324 \pm 1.1 \times 10^{-5}$, occurring when both $k$ and $l$ follow beta distributions. More generally, changing the distribution of $l$ appears to have a greater effect than changing the distribution of $k$: the estimates are consistently relatively high when $l$ is log-normally distributed, whereas no distribution of $k$ produces consistently high estimates across all distributions of $l$. Although the estimates vary substantially between combinations, they remain small in absolute terms, with the largest corresponding to $1.6\%$.

Next, let \(X_1, X_2, \ldots, X_5\) be independent random variables with \(X_i\sim\operatorname{Binomial}(9,0.5)\), and define \(h_i=1+X_i\).
Table~\ref{tab:kl-binomial_hi} shows a similar pattern to Table~\ref{tab:kl-uniform_hi}. For every distribution of $k$, the highest estimated probability occurs when $l$ follows a log-normal distribution. The largest estimate is $0.007744 \pm 5.4 \times 10^{-5}$, occurring when $k$ follows a beta distribution, while the lowest estimate is $0.0000274 \pm 3.2 \times 10^{-6}$, occurring when $k$ follows a beta distribution and $l$ follows an Erlang-16 distribution. Once again, the distribution of $l$ appears to have a greater effect than the distribution of $k$. All estimated probabilities in this 
table are below $0.8\%$.

\begin{table}[H]
\centering
\renewcommand{\arraystretch}{1.6}
\resizebox{\textwidth}{!}{%
\begin{tabular}{|c||c|c|c|c|c|c|c|}
\hline
\raisebox{-3mm}{
\textbf{Distribution}} 
& \multicolumn{7}{c|}{\textbf{Distribution of $l$:}} \\
\cline{2-8}
\textbf{of $k$:} & \textbf{Uniform} & \textbf{Exponential} & \textbf{Erlang-2} & \textbf{Erlang-4} & \textbf{Erlang-16} & \textbf{Log-normal} & \textbf{Beta} \\
\hline\hline
\textbf{Uniform} & 0.0019 & 0.0031 & 0.0016 & 0.00077 & 0.00025 & 0.0067 & 0.00023 \\
\hline
\textbf{Exponential} & 0.0018 & 0.0029 & 0.0017 & 0.00097 & 0.00053 & 0.0058 & 0.00049 \\
\hline
\textbf{Erlang-2} & 0.0019 & 0.0031 & 0.0016 & 0.00082 & 0.00034 & 0.0066 & 0.00032 \\
\hline
\textbf{Erlang-4} & 0.0019 & 0.0032 & 0.0015 & 0.00067 & 0.00018 & 0.0072 & 0.00017 \\
\hline
\textbf{Erlang-16} & 0.0018 & 0.0032 & 0.0014 & 0.00048 & 0.000037 & 0.0077 & 0.000040 \\
\hline
\textbf{Log-normal} & 0.0014 & 0.0019 & 0.0013 & 0.00097 & 0.00073 & 0.0035 & 0.00072 \\
\hline
\textbf{Beta} & 0.0018 & 0.0033 & 0.0014 & 0.00046 & 0.000027 & 0.0077 & 0.000035 \\
\hline
\end{tabular}%
}
\caption{Estimated probability of Braess' paradox under different distributions of $k$ and $l$, with the hierarchy parameters independently following a $1+\mathrm{Binomial}(9,0.5)$ distribution.}
\label{tab:kl-binomial_hi}
\end{table}

Finally, let the hierarchy parameters be independent and geometrically distributed with probability \(1/5.5\).
Table~\ref{tab:kl-geometric_hi} again shows relatively high estimated probabilities when $l$ follows a log-normal distribution. The largest estimate is $0.009390 \pm 6.0 \times 10^{-5}$, occurring when $k$ follows a beta distribution, while the lowest estimate is $0.000624 \pm 1.5 \times 10^{-5}$, occurring when both $k$ and $l$ follow beta distributions. A notable exception to the wider variation across distributions of $l$ occurs when $k$ is log-normally distributed: the estimates then remain between $0.0070$ and $0.0073$ for all distributions of $l$. All estimated probabilities are below $1\%$.

\begin{table}[H]
\centering
\renewcommand{\arraystretch}{1.6}
\resizebox{\textwidth}{!}{%
\begin{tabular}{|c||c|c|c|c|c|c|c|}
\hline
\raisebox{-3mm}{
\textbf{Distribution}} 
& \multicolumn{7}{c|}{\textbf{Distribution of $l$:}} \\
\cline{2-8}
\textbf{of $k$:} & \textbf{Uniform} & \textbf{Exponential} & \textbf{Erlang-2} & \textbf{Erlang-4} & \textbf{Erlang-16} & \textbf{Log-normal} & \textbf{Beta} \\
\hline\hline
\textbf{Uniform} & 0.0051 & 0.0063 & 0.0050 & 0.0041 & 0.0034 & 0.0091 & 0.0033 \\
\hline
\textbf{Exponential} & 0.0061 & 0.0068 & 0.0059 & 0.0053 & 0.0050 & 0.0089 & 0.0050 \\
\hline
\textbf{Erlang-2} & 0.0053 & 0.0064 & 0.0050 & 0.0042 & 0.0036 & 0.0093 & 0.0036 \\
\hline
\textbf{Erlang-4} & 0.0044 & 0.0059 & 0.0042 & 0.0031 & 0.0022 & 0.0094 & 0.0021 \\
\hline
\textbf{Erlang-16} & 0.0035 & 0.0054 & 0.0034 & 0.0021 & 0.00083 & 0.0094 & 0.00073 \\
\hline
\textbf{Log-normal} & 0.0073 & 0.0072 & 0.0070 & 0.0070 & 0.0071 & 0.0073 & 0.0071 \\
\hline
\textbf{Beta} & 0.0034 & 0.0054 & 0.0034 & 0.0020 & 0.00073 & 0.0094 & 0.00062 \\
\hline
\end{tabular}%
}
\caption{Estimated probability of Braess' paradox under different distributions of $k$ and $l$, with the hierarchy parameters independently following a $\mathrm{Geometric}(1/5.5)$ distribution.}
\label{tab:kl-geometric_hi}
\end{table}

Across the three tables, a log-normal distribution of $l$ is consistently associated with the highest estimated probabilities. The distribution of $l$ generally appears to have a greater influence than that of $k$, demonstrating that the roles of the two parameters are not symmetric. The shifted binomial distribution of the hierarchy parameters produces consistently lower estimates than the uniform and geometric distributions for the combinations considered. Very low probabilities were observed when both $k$ and $l$ followed a beta distribution. The highest probability estimate across Tables~\ref{tab:kl-uniform_hi}--\ref{tab:kl-geometric_hi} is $1.6\%$, 
obtained with uniformly distributed hierarchy parameters.

\section{Conclusions}

We have shown that Braess' paradox can occur in hierarchical network configurations that could potentially be used to model the British road system. We derived necessary and sufficient conditions for its occurrence at some total flow \(Q>0\), identified cases in which it is impossible, and gave sufficient conditions and explicit examples in which it occurs. We also obtained probability bounds for randomly assigned hierarchy levels when $k$ is fixed and $l$ may vary.

We used Monte Carlo simulations to estimate the probability of the paradox occurring under different distributions of the network parameters. The results suggest that, although the probability is generally low, it is sensitive to the choice of parameter distributions. Simulations using beta distributions for $k$ and $l$ reveal a consistent association between greater positive skewness and a higher probability of Braess' paradox. 

Our results support the general conjecture stated by \cite{VZverovich2021} that the probability of Braess' paradox occurring is typically small. 
Under the distributions considered, all Monte Carlo estimates were below 3\%, with most falling below 1.6\%.  
In particular, when $k$ and $l$ follow a shifted Erlang-2 distribution, the probability of Braess' paradox occurring in the hierarchical network ranges from 0.16\% to 0.50\%, depending on the distributions of the hierarchy parameters $h_i$. In contrast, the corresponding probability is 6\% in the classical network, where the parameters are modelled using an Erlang-2 distribution \citep{VZverovich2021}. 

Moreover, these results potentially suggest that, from the perspective of congestion minimisation, hierarchical networks may be more suitable for modelling real-world transportation systems than standard networks, for which the probability of Braess' paradox occurring can be as high as 14\% (see \cite{VZverovich2021}).

\bibliography{bibliography}

\appendix

\section{Proof of Corollary~\ref{cor:Probability_dist_hi}}
\label{app:Probability_dist_hi_proof}

\begin{proof}
Let \(\mathcal{B}_k\) denote the event that Braess' paradox occurs in
\(V/V^+\) for some \(Q>0\) with the given hierarchy parameters and for some
$l>1$. Thus,
\[
P_k=\mathbb{P}(\mathcal{B}_k).
\]
Recall that
\[
A:=k^{h_2}-k^{h_3}-k^{h_5},
\qquad
B:=k^{h_4}-k^{h_1}-k^{h_3}.
\]

We first prove the lower bound. Define
\[
E_A:=
\left\{
\phi_A<h_2-\log_k 2
\text{ and }
h_1-h_2>h_4-h_5>0
\right\}
\]
and
\[
E_B:=
\left\{
\phi_B<h_4-\log_k 2
\text{ and }
h_5-h_4>h_2-h_1>0
\right\}.
\]
Suppose that \(E_A\) occurs. Since \(h_3,h_5\leq\phi_A\), we have
\[
A
=k^{h_2}-k^{h_3}-k^{h_5}
\geq k^{h_2}-2k^{\phi_A}
>0,
\]
where the final inequality follows from
\(\phi_A<h_2-\log_k 2\). 
Let us prove \eqref{eq:HBP1}:
\[
h_1+h_5-(h_2+h_4)
=(h_1-h_2)-(h_4-h_5)>0.
\]
Hence, \(h_1+h_5>h_2+h_4\), while
$
h_1>h_2
$
and
$
h_4>h_5.
$
Therefore, condition~\eqref{eq:HBP1} and Case~\(1\) of
Theorem~\ref{thm:Braess_existence} are satisfied. Consequently,
$
E_A\subseteq\mathcal{B}_k.
$

Similarly, suppose that \(E_B\) occurs. Since \(h_1,h_3\leq\phi_B\),
\[
B
=k^{h_4}-k^{h_1}-k^{h_3}
\geq k^{h_4}-2k^{\phi_B}
>0.
\]
Further, let us prove \eqref{eq:HBP1}:
\[
h_1+h_5-(h_2+h_4)
=(h_5-h_4)-(h_2-h_1)>0
\]
and also $h_1<h_2$
and $h_4<h_5.$
Thus, condition~\eqref{eq:HBP1} and Case~\(2\) of
Theorem~\ref{thm:Braess_existence} hold, so
$
E_B\subseteq\mathcal{B}_k.
$

The events \(E_A\) and \(E_B\) are disjoint, since \(E_A\) requires
\(h_1>h_2\) and \(h_4>h_5\), whereas \(E_B\) requires
\(h_1<h_2\) and \(h_4<h_5\). Therefore,
\[
\mathbb{P}(E_A)+\mathbb{P}(E_B)
=\mathbb{P}(E_A\cup E_B)
\leq\mathbb{P}(\mathcal{B}_k)
=P_k.
\]

For the upper bound, define
\[
F:=
\left\{
h_1+h_5 \leq h_2+h_4
\right\}
\cup
\left\{
h_2-\log_k 2 \leq \psi_A
\ \text{and}\
h_4-\log_k 2 \leq \psi_B
\right\}.
\]
If \(h_1+h_5 \leq h_2+h_4\), then the necessary
condition~\eqref{eq:HBP1} fails, and hence Braess' paradox cannot occur for any $l>1$.
Alternatively, suppose that
\[
h_2-\log_k 2 \leq \psi_A=\min\{h_3,h_5\},
\qquad
h_4-\log_k 2 \leq \psi_B=\min\{h_1,h_3\}.
\]
Then
\[
A
=k^{h_2}-k^{h_3}-k^{h_5}
\leq k^{h_2}-2k^{\psi_A}
\leq 0,
\]
and, similarly,
\[
B
=k^{h_4}-k^{h_1}-k^{h_3}
\leq k^{h_4}-2k^{\psi_B}
\leq 0.
\]
Consequently, for every $l>1$,
\[
l^{h_1}A+l^{h_2}B\leq 0
\qquad\text{and}\qquad
l^{h_4}A+l^{h_5}B\leq 0.
\]
Therefore, conditions~\eqref{eq:HBP2} and~\eqref{eq:HBP3} cannot
be satisfied, so Braess' paradox cannot occur for any $l>1$.

It follows that
\[
F\subseteq\mathcal{B}_k^{\,c}.
\]
Hence,
\[
\mathbb{P}(F)\leq\mathbb{P}(\mathcal{B}_k^{\,c})=1-P_k
\]
and, therefore,
\[
P_k
\leq
1-\mathbb{P}\!\left[
h_1+h_5\leq h_2+h_4
\ \text{or}\
\bigl(
h_2-\log_k 2\leq\psi_A
\ \text{and}\
h_4-\log_k 2\leq\psi_B
\bigr)
\right].
\]
This completes the proof.
\end{proof}

\end{document}